\documentclass[12pt,leqno]{article}
\usepackage{amssymb,amsmath,amsfonts,amsbsy, xspace, latexsym}
\newtheorem{theorem}{Theorem}
\newtheorem{lemma}{Lemma}

\newtheorem{definition}{Definition}
\newtheorem{corollary}{Corollary}
\newtheorem{notation}{Notation}
\newtheorem{remark}{Remark}
\newtheorem{example}{Example}

\newcommand{\N}{\ensuremath{\,\mathbb{N}}}
\newcommand{\Z}{\ensuremath{\,\mathbb{Z}}}

\newcommand{\R}{\ensuremath{\,\mathbb{R}}}
\newcommand{\C}{\ensuremath{\,\mathbb{C}}}

\newcommand{\ifff}{\ensuremath{\mbox{\; if{f}\; }}}

\newcommand{\Proof}{\noindent \textbf{Proof:\;}}

\newcommand{\bra}{\ensuremath{\large\langle}}

\newcommand{\ket}{\ensuremath{\large\rangle}}

\newcommand{\dom}{\ensuremath{\rm{dom}}}

\newcommand{\ran}{\ensuremath{\rm{ran}}}

\newcommand{\Span}{\ensuremath{\rm{Span}}}

\newcommand{\qst}{\ensuremath{\,\widehat{\rm{st}}}}

\begin{document}
{\small
%___________________________________________________________________________

\title {Existence and Uniqueness of $v$-Asymptotic Expansions and 
Colombeau's Generalized Numbers}

\author{Todor D. Todorov \\
			Mathematics Department\\		
			California Polytechnic State University\\
			San Luis Obispo, California 93407, USA\\
(ttodorov@calpoly.edu)}

\date{}
\maketitle	

\begin{abstract}
We define a type of generalized asymptotic series called $v$-asymptotic. We show
that every function with moderate growth at infinity has a $v$-asymptotic expansion. We also describe
the set of $v$-asymptotic series, where a given function with moderate growth has a unique
$v$-asymptotic expansion. As an application to random matrix theory we calculate the
coefficients and establish the uniqueness of the $v$-asymptotic expansion of an integral 
with a large
parameter. As another application (with significance in the non-linear theory of generalized
functions) we show that every Colombeau's generalized number has a
$v$-asymptotic expansion. A similar result follows for Colombeau's generalized functions, in particular,
for all Schwartz distributions. 
\end{abstract}
{\small Mathematics Subject Classification: 30B10, 34E05, 35D05, 41A60, 40A30, 46F30.
Key words: asymptotic expansion, valuation, ultrametric space, Colombeau generalized functions,
random matrix theory.}
\section{Introduction}\label{S: Introduction}

 Our framework is the ring $\mathcal{M}$ of
functions with moderate growth at infinity, supplied with a pseudovaluation $v: \mathcal{M}\to
\R\cup\{\infty\}$ and pseudometric \newline$d_v: \mathcal{M}^2\to\R$ (Section~\ref{S:
Functions of Moderate Growth}). We denote by
$\mathcal{F}_v$ the ring of functions with the non-negative pseudovaluation and by $\mathcal{C}_v$ the
set of functions with zero pseudovaluation. Here are some typical examples:
$r+1/x^n\in\mathcal{F}_v,\; r\in\R,\; n=0, 1, 2,\dots$ and\; $\sin{x}, \cos{x}, \ln{x}\in\mathcal{C}_v$.

	In Section~\ref{S: Pseudostandard Part} we prove the existence of a linear
homomorphism (a linear operator) $\qst$, called the ``pseudostandard part mapping,'' from
$\mathcal{F}_v$ into $\mathcal{F}_v$ with range within $\mathcal{C}_v\cup\{0\}$. The mapping $\qst$ is an
extension of the limit $\lim_{x\to\infty}$, considered as a functional (part~(iii) of Theorem~\ref{T:
Properties}).

	In Section~\ref{S: v-Asymptotic Series} we define the concept of $v$-asymptotic series; these are
series of the form $\sum_{n=0}^\infty\; \frac{\varphi_n(x)}{x^{r_n}}$, where $(r_n)$ is a strictly
increasing (bounded or unbounded) sequence in
\R\, and the set $\{\varphi_n(x)\mid n=0,1,\dots\}$ is a linearly $v$-independent subset of
$\mathcal{C}_v$ (Section~\ref{S: Linearly v-Independent Sets}). In that sense, the
functions in $\mathcal{C}_v$ play the role of  ``our generalized numbers.'' There are many examples in
this section: some series are both $v$-asymptotic and asymptotic (in the usual sense), some series are
$v$-asymptotic, but not asymptotic and others are asymptotic, but not $v$-asymptotic.

	The main result is in Theorem~\ref{T: Existence and Uniqueness} (Section~\ref{S: Existence
and Uniqueness Result}): Every function with moderate growth at infinity has a
$v$-asymptotic expansion. We also describe
the sets of $v$-asymptotic series, where every function in $\mathcal{M}$ has a unique
$v$-asymptotic expansion. We
shall try to explain the idea of the proof by modifying a familiar theorem in asymptotic analysis: Let
$f(x)$ be a function, $(r_n)$ be a strictly increasing sequence in
\R\, and $(c_n)$ be a sequence in \C. The asymptotic expansion $f(x)\thicksim
\sum_{n=0}^\infty\; \frac{c_n}{x^{r_n}},\, x\to\infty$,  holds {\em if{f}}
\begin{equation}\label{E: Limits}
c_n=\lim_{x\to\infty}\left(x^{r_n}
\left[f(x)-\sum_{k=0}^{n-1}\frac{c_k}{x^{r_k}}\right]\right),
\end{equation}
for $n=0, 1, 2,\dots$. This result can be rephrased as follows: If all limits (for $n=0, 1, 2,\dots$) on
the RHS of (\ref{E: Limits}) exist, then $f(x)$ has a unique asymptotic expansion of the form
$\sum_{n=0}^\infty\;\frac{c_n}{x^{r_n}}$. It is clear that most of the functions in $\mathcal{M}$ do not
have asymptotic expansions of this form. Actually, there are very few general existence results
in asymptotic analysis
(N.~G.~De~Bruijn~\cite{deBruijn}, R.~Estrada~and~R.~P.~Kanwal~\cite{EstradaKanwal}). This is,
perhaps, the reason why many consider the asymptotic analysis as a ``kind of art'' or ``collection
of different methods,'' rather than a mature mathematical theory. Our article is an attempt to
improve this situation. Here is the way we modify the above asymptotic theorem: 

	(a) Instead of series of the form  $\sum_{n=0}^\infty\; \frac{c_n}{x^{r_n}}$, we consider
$v$-asymptotic series
$\sum_{n=0}^\infty\; \frac{\varphi_n(x)}{x^{r_n}}$, where $v(\varphi_n)=0$ and $v$ is the
pseudovaluation mentioned earlier.

	(b) We replace $\lim_{x\to\infty}$ (considered as a functional) in the counterpart of (\ref{E: Limits})
by the linear operator $\qst$ mentioned earlier. The advantage of $\qst$ over $\lim_{x\to\infty}$ is that
$\qst$ is defined on the whole $\mathcal{F}_v$ (including functions such as $\sin{x}, \cos{x}, \ln{x},
\ln{(\ln{x})}$, etc.). In contrast,
$\lim_{x\to\infty}$ is defined only on a proper subset of $\mathcal{F}_v$ (e.g. $\lim_{x\to\infty}\sin{x}$
and $\lim_{x\to\infty}\ln{x}$ do not exist). 

	(c) Instead of a (fixed) sequence $(r_n)$ (which corresponds to a choice of a fixed asymptotic scale), we
allow different sequences $(r_n)$ for the different functions $f(x)$. 

	Two applications are presented in the article: At the end of Section~\ref{S: Existence and Uniqueness
Result} we prove the uniqueness of the $v$-asymptotic expansion of an integral with a large parameter
with origin in random matrix theory. In Section~\ref{S: Colombeau's Generalized Numbers} we show that
every Colombeau's generalized number has a $v$-asymptotic expansion. A similar result follows for
Colombeau's generalized functions, in particular, for all Schwartz distributions. This result has
importance in the nonlinear theory of generalized functions (J.F.
Colombeau~\cite{jfCol84}-\cite{jfCol85}) and its applications in ordinary and partial
differential equations (M. Oberguggenberger~\cite{mOber92}). 

	This article has many features in common with
another article by T.D. Todorov~and R. Wolf~\cite{TodWolf} on A. Robinson's asymptotic numbers in
the framework of the non-standard asymptotic analysis (A.~H.~Lightstone~and~A.
Robinson~\cite{LiRob}).

	In what follows $\mathbb{N}=\{1, 2, 3, \dots, \}$ denotes the set of natural numbers and
$\mathbb{R}$ and $\mathbb{C}$ denote the fields of the real and complex numbers, respectively.
%----------------------------
\section{Our Framework: Functions of Moderate Growth at Infinity}\label{S: Functions of Moderate Growth}

We define and study the properties of a pseudovaluation and the corresponding pseudometric in the ring of
complex valued functions with moderate growth at infinity. Most of the results are
elementary and presented without proofs.

\begin{definition}[Functions of Moderate Growth]\label{D: Functions of Moderate Growth} 

	{\bf (i)} We denote by $G_\infty(\mathbb{R}, \mathbb{C})$ the ring (under the usual pointwise addition
and multiplication) of all complex valued functions $f$ defined on
$\dom{(f)}\subseteq\mathbb{R}$ such that $(n, \infty)\subseteq \dom{(f)}$ for some
$n\in \mathbb{N}$. We shall treat the elements of $G_\infty(\mathbb{R}, \mathbb{C})$ as germs at $\infty$,
i.e. two functions will be identified if they have the same values for all sufficiently large $x$. We define
a partial order in $G_\infty(\mathbb{R}, \mathbb{C})$ by
$f>0$ if $f(x) > 0$ for all sufficiently large\, $x$.

  {\bf (ii)} The functions in 
\[
\mathcal{M} = \{\, f\in G_\infty(\mathbb{R}, \mathbb{C})\,\mid\, f(x) = O(x^n),\mbox{\; as\;}
x\to\infty,\mbox{\; for some\;\,} n\in\N\,\}
\]
are called {\bf functions with moderate growth at $\infty$} (or, {\bf moderate functions}, for short) and
those in
\[
\mathcal{N} = \{\, f\in G_\infty(\mathbb{R}, \mathbb{C})\,\mid\, f(x) = O(1/x^n),\mbox{\;as\;}
x\to\infty, \mbox{\; for all\quad} n\in\N\,\}
\]
are called {\bf null-functions}.

	{\bf (iii)} We define a {\bf pseudovaluation}\; $v: \mathcal{M}\to\R\cup\{\infty\}$ by
\begin{equation}\label{E: PsuedoValuation}
			v(f) = \sup\{\, r\in\R\, \mid\, f(x) = O(1/x^r) \mbox{\; as\quad} x\to\infty\,\}.
\end{equation}
The functions in the sets: 
\begin{align} 
	&\mathcal{I}_v = \{\, f\in\mathcal{M}\, \mid\, v(f) > 0\, \},  \label{E: v-infinitesimals} \\
	&\mathcal{F}_v = \{\, f\in\mathcal{M}\, \mid\, v(f) \geq 0\,\},  \label{E: v-finite}\\
	&\mathcal{C}_v = \{\, f\in\mathcal{M}\, \mid\, v(f) = 0\,\},  \label{E: v-constants}
\end{align}
are called {\bf $v$-infinitesimal, $v$-finite} and {\bf $v$-constants}, respectively. The functions in
$\mathcal{M}\setminus\mathcal{F}_v$ will be called {\bf $v$-infinitely large}.

	{\bf (iv)} If $S\subseteq G_\infty(\mathbb{R}, \mathbb{C})$, then the set:
\begin{equation}\label{E: Monad}
\mu_v(S)=\{\varphi+d\varphi\mid \varphi\in S,\; d\varphi\in \mathcal{I}_v\, \}
\end{equation}
is called the {\bf
$v$-monad} of $S$.  In the particular case $S=\{f\}$, we shall write $\mu_v(f)$ instead of the more
precise $\mu_v(\{f\})$.
\end{definition}

\begin{example}: $x^r,\;  e^{ix},\; \sin{x},\; \cos{x},\; \ln{x},\; \ln({\ln{x})}$ are all in
$\mathcal{M}$ (where
$r
\in
\R$). Also, $0$ and $e^{-x}$ are in $\mathcal{N}$. In contrast, 
$e^x\notin\mathcal{M}$. Notice that $\mathcal{M}$ is subring of $G_\infty(\mathbb{R}, \mathbb{C})$ with
zero divisors. In particular, the functions  $x^r,\; e^{ix},\; \ln{x},\; \ln({\ln{x})}$ are all
multiplication invertible in
$\mathcal{M}$, while $\sin{x}, \cos{x}$ are not. 
\end{example} 

\begin{example}: If $r \in \R$, then $v(1/x^r) = r$. More generally, $v(P) = -\deg(P)$ for any
polynomial $P\in\C[x]$, where $\deg(P)$ denotes the degree of $P$. We have
$v(c)=0$ for all
$c\in\C,\; c\neq 0$. In other words,
$\C\setminus\{0\}\subset\mathcal{C}_v$. Also, $\sin{x},\;
\cos{x},\; e^{ix},\; \ln{x},\; \ln^n{x}$ and $\ln({\ln{x})}$ are all in $\mathcal{C}_v$, because
$v(e^{ix}) = v(\sin{x}) = v(\cos{x}) = v(\ln{x}) = v(\ln^n{x})= v(\ln({\ln{x})}) = 0$ ($n\in\Z$). 
  Finally, $v(0) =
v(e^{-x}) =
\infty$.
\end{example} 
	
	The results of the next three lemmas follow immediately from the definition of the pseudovaluation:

\begin{lemma}[A Characterization of $v$]\label{L: Characterization of v} Let $f\in
G_\infty(\mathbb{R}, \mathbb{C})$ and $r\in\R$. Then the following are equivalent:
	
	{\bf (i)}\;   $v(f) = r$.

	{\bf (ii)}	$\frac{1}{ x^r\sqrt[n]{x}} \leq |f(x)|\leq \frac{\sqrt[n]{x}}{ x^r}$ for all $n\in \N$.

	{\bf (iii)} $f(x) = \frac{\varphi(x)}{x^r}$ for some $\varphi \in \mathcal{C}_v$.

\end{lemma}
\begin{lemma}[A Characterization of $\mathcal{N}$] Let $f\in G_\infty(\mathbb{R}, \mathbb{C})$. The
following are equivalent:
	
	{\bf (i)}\quad   $f \in \mathcal{N}$.

	{\bf (ii)}\,  $v(f) = \infty$.

	{\bf (iii)}  $\lim_{x\to\infty} x^nf(x) = 0$ for all $n\in\N$.

\end{lemma}
	
\begin{lemma}[Connection with O's Symbols]\label{L: Connection with O's Symbols}
	
	{\bf (i)}\quad   $\mathcal{I}_v \subset o(1)\subset \mathcal{F}_v$ in the sense that $f \in
\mathcal{I}_v$ implies $f(x) = o(1)$, as $x\to\infty$, which implies $f\in\mathcal{F}_v$.

	{\bf (ii)}\,  $O(1)\setminus o(1)\subset \mathcal{C}_v$ in the sense that if $f(x) = O(1)$ and 
$f(x) \neq o(1)$, as $x\to\infty$, then  $f \in\mathcal{C}_v$. More
generally,  $f\in O(1/x^r)\setminus o(1/x^r)$ implies
$v(f) = r$.

	{\bf (iii)} The functions in $G_\infty(\mathbb{R}, \mathbb{C})\setminus\mathcal{F}_v$ are unbounded
on every interval of the form $(n, \infty),\, n \in \N$.

\end{lemma}

	Part (i) shows that the $v$-infinitesimals are proper infinitesimals (with respect to the order relation
in $\mathcal{M}$). Part (ii) shows that the finite but non-infinitesimal functions (with respect to the
order in
$\mathcal{M}$) are $v$-finite. Part (iii)  shows that if $f$ is outside $\mathcal{F}_v$ and
if $|f|$ happens to be in order with all $n\in\N$, then $f$ is an infinitely large element of
$G_\infty(\mathbb{R}, \mathbb{C})$. We should notice that $\mathcal{C}_v$ contains also infinitely large
functions (along with finite and infinitesimal). For example,  $\ln{x}$ is infinitely large in the
sense that $n <\ln{x}$ for all $n\in\N$ and $1/\ln{x}$ is a positive infinitesimal in the sense that
$0<1/\ln{x}<1/n$ for all $n\in\N$.

\begin{theorem}[Properties of $v$]\label{T: Properties of v} The function $v$
is a pseudovaluation on
$\mathcal{M}$ in the sense that:
	
	{\bf (i)}\quad $v(f) = \infty$ \ifff $f\in\mathcal{N}$.

	{\bf (ii)}\; $v(f g) \geq v(f) + v(g)$. In particular, we have  $v(c g)= v(f)$ for all $c\in\C,\; c\not=0$ and also  $v(f/x^r)= v(f)+r$ for all $r\in\R$.

	{\bf (iii)}\;\, $v(f \pm g) \geq \min\{v(f), v(g)\}$. Moreover, $v(f)\not=v(g)$ implies\newline $v(f
\pm g) = \min\{v(f), v(g)\}$.

	In addition, $v$ has the following properties: 

	{\bf (iv)} $v(-f) = v(f) = v(|f|)$;
	
	{\bf (v)}\, $v(1/f) = -v(f)$ whenever $f$ is invertible in $\mathcal{M}$;

	{\bf (vi)}\, $|f|<|g|$ implies $v(f) \geq v(g)$.
\end{theorem}
\Proof The above properties follow directly from the properties of the ``$\sup$'' in the definition of $v$
and we leave the verification to the reader. $\blacktriangle$
\begin{theorem}[Rings and Ideals]\label{T: Rings and Ideals}

	{\bf (i)} $\mathcal{M}$ is convex subring (with zero divisors) of $G_\infty(\mathbb{R}, \mathbb{C})$,
where ``convex'' means that  if $f\in G_\infty(\mathbb{R}, \mathbb{C})$ and $g
\in\mathcal{M}$, then $|f| \leq |g|$ implies $f\in \mathcal{M}$. Also  $\mathcal{N}$ is a convex ideal in
$\mathcal{M}$.

	{\bf (ii)} We have 	$\mathcal{F}_v =\mathcal{C}_v \cup \mathcal{I}_v$\; and\; $\mathcal{C}_v \cap
\mathcal{I}_v =\varnothing$. Besides,
$\mathcal{F}_v$ is convex subring of $\mathcal{M}$ and $\mathcal{I}_v$ is a convex ideal in
$\mathcal{F}_v$.

	{\bf (iii)} The ring $\mathcal{F}_v$ contains a copy of the field of complex numbers
\C\, presented by the constant-functions. Also we have  $\C\subset\mathcal{C}_v\cup\{0\}$.
\end{theorem}

\Proof (i) follows immediately from the definitions of $\mathcal{M}$ and $\mathcal{N}$.

	(ii) Suppose $f, g\in\mathcal{F}_v$, i.e. $v(f)\geq 0, v(g)\geq 0$.
It follows $v(f+g)\geq\min\{v(f),\; v(g)\}\geq 0$ and $v(fg)\geq
v(f)+v(g)\geq 0$, by Theorem~\ref{T: Properties of v}, i.e. $f+g,\; fg\in\mathcal{F}_v$. Similarly it
follows that $\mathcal{I}_v$ is closed under the addition and multiplication. Suppose also
that $h\in\mathcal{I}_v$, i.e. $v(h)>0$. We have
$v(fh)\geq v(f)+v(h)>0$, i.e. $fh\in\mathcal{I}_v$, as required. 

	(iii) follows from the fact that $v(c)=0$ for all $c\in\C,\, c\not=0$.
$\blacktriangle$
\begin{definition}[Pseudometric]\label{D: Pseudometric}  We define a pseudometric\;
$d_v: \mathcal{M}^2\to\R$ by $d_v(f, g) = e^{-{v(f-g)}}$ (with the understanding that
$e^{-\infty}= 0$). We denote by $(\mathcal{M}, d_v)$ the corresponding pseudometric space. 
\end{definition}

\begin{theorem}[Properties of the Pseudometric]\label{T: Properties of the Pseudo-Metric} For
any $f, g, h \in \mathcal{M}$:
	
	{\bf (i)}\;\; $d_v(f, g) = 0$ \ifff  $f-g\in\mathcal{N}$.

	{\bf (ii)}\; $d_v(f, g) = d_v(g, f)$.

	{\bf (iii)} $d_v(f, g) \leq \max\{d_v(f, h), d_v(h, g)\}\leq d_v(f, h) + d_v(h, g)$. 

	{\bf (iv)}  $\mathcal{I}_v, \mathcal{F}_v$ and $\mathcal{C}_v$ coincide with the open unit ball, the closed unit ball and the unit sphere in
$(\mathcal{M}, d_v)$, respectively.
\end{theorem}
\Proof A direct consequence of the property of the pseudovaluation $v$ (Theorem~\ref{T: Properties of
v}).
	We should mention that $d_v(f, g) \leq \max\{d_v(f, h), d_v(h, g)\}$ is called {\bf ultrametric
inequality}.
\begin{theorem}[Convergence in $\mathcal{M}$]\label{T: Convergence in M} Let $(f_n)$ be
a sequence in $\mathcal{M}$.

	{\bf (i)} If $f\in\mathcal{N}$ (in particular, $f=0$), then $\lim_{n\to\infty}f_n = f$ in $(\mathcal{M},
d_v)$ \ifff $\lim_{n\to\infty}v(f_n) = \infty$ in
$\R\cup\{\infty\}$.

	{\bf (ii)} If $f\notin\mathcal{N}$, then $\lim_{n\to\infty}f_n = f$ in $(\mathcal{M}, d_v)$
implies that $v(f_n) = v(f)$ for all sufficiently large $n$.
\end{theorem}
\Proof 

	(i) $\lim_{n\to\infty}f_n = f$ in $(\mathcal{M}, d_v)$ \ifff $\lim_{n\to\infty}d(f_n, f)=0$ in
\R\, \ifff\newline $\lim_{n\to\infty}e^{-v(f_n-f)}=0$ in $\R$ \ifff
$\lim_{n\to\infty}v(f_n-f)=\infty$ in $\R\cup\{\infty\}$
\ifff\newline $\lim_{n\to\infty}v(f_n)=\infty$, as required, since $v(f_n-f)=v(f_n)$ for all sufficiently
large $n$, by part~(iii) of Theorem~\ref{T: Properties of v}. 

	(ii) We have $\lim_{n\to\infty}v(f_n-f)=\infty$, by (i). Suppose (on the contrary)
that there exists a unbounded sequence
$(\nu_n)$ in $\N$ such that
$v(f_{\nu_n})\not=v(f)$ for all $n$. It follows $v(f_{\nu_n}-f)=\min\{v(f_{\nu_n}),v(f)\}$ for all $n$, by
part~(iii) of Theorem~\ref{T: Properties of v}, which implies $\lim_{n\to\infty}v(f_{\nu_n})=\infty$.
Finally, it follows $v(f)=\infty$, a contradiction.
$\blacktriangle$
\begin{corollary}\label{C: For LC} Let $(r_n)$ be a sequence in $\R\cup\{\infty\}$ such
that $\lim_{n\to\infty}r_n=\infty$. Let $(f_n)$ be a sequence in $\mathcal{M}$ such that
$f_n=o(1/x^{r_n})$ or $f_n=O(1/x^{r_n})$ as $x\to\infty$ for all $n$. Then
$\lim_{n\to\infty}f_n=0$ in $(\mathcal{M}, d_v)$.
\end{corollary}
\Proof The result follows, by part (i) of the above theorem, since in both cases $v(f_n)\geq r_n$ 
for all sufficiently large $n$.
$\blacktriangle$
\begin{example} We have $\lim_{n\to\infty}\frac{1}{x^n}=0$ in $(\mathcal{M}, d_v)$ since
$v(\frac{1}{x^n})=n\to\infty$ as $n\to\infty$.
\end{example}
%--------------------------------------
\section{Linearly $v$-Independent Sets}\label{S: Linearly v-Independent Sets}

 We should notice that the set of the $v$-constants $\mathcal{C}_v$ is not closed
under addition and multiplication in $\mathcal{M}$. Our next goal 
is to look for vector spaces inside $\mathcal{C}_v\cup\{0\}$.

	If $\Phi\subset \mathcal{M}$, we denote by $\Span(\Phi)$
the span of $\Phi$ over $\mathbb{C}$ (within $\mathcal{M}$). Recall that $\Span(\Phi)$ is the smallest
vector subspace of
$\mathcal{M}$ containing $\Phi$.
\begin{definition}[$v$-Independence]\label{D: v-Independence} We say that a set $\Phi\subset
\mathcal{M}$ is {\bf linearly $v$-independent over \C} if\, 
$\Span(\Phi)\subset\mathcal{C}_v
\cup\{0\}$.
\end{definition}

	Let $\Phi$ be an linearly $v$-independent set over \C. It is clear that
$\Phi\subset\mathcal{C}_v\cup\{0\}$. Also, every subset $\Psi$ of $\Phi$ is also linearly
$v$-independent over \C.

 Here is another (trivial) reformulation of the above definition:

\begin{lemma}\label{L: Trivial} Let $\Phi \subset \mathcal{M}$. Then the following are equivalent:

	{\bf (i)}\quad  $\Phi$ is linearly $v$-independent over \C.
	
	{\bf (ii)}\quad $v(f) = 0$ for all $f\in \Span(\Phi),\, f\not=0$.

	{\bf (iii)}\, There exists a vector subspace $K$ of $\mathcal{M}$ such that
$\Phi\subseteq K\subset\mathcal{C}_v\cup\{0\}$.
\end{lemma} 
\begin{example} \label{Ex: v-independent} The set of the complex numbers $\mathbb{C}$ is linearly
$v$-independent (in a trivial way). If
$f\in\mathcal{C}_v\cup\{0\}$, then the set $\{f\}$ is also linearly $v$-independent over \C. In
particular, the set $\{0\}$ is linearly $v$-independent over \C. 
\end{example}
\begin{lemma} \label{L: v-Independent Sets} Each of the following sets is linearly $v$-independent:
$\Phi_1 = \{\sin{x},\; \cos{x}\, \}$,\; $\Phi_2 = \{\ln{x},\; \ln({\ln{x}}),\; 
\ln(\ln({\ln{x}})),\dots\}$,\; $\Phi_3 = \{\ln{x},\; \ln^2{x},\; \ln^3{x},\; \ln^4{x},
\dots\}$,\; $\Phi_4 = \{e^{\pm i\pi x},\quad e^{\pm i\pi x}\, \ln{x},\quad e^{\pm i\pi x}\,
\ln^2{x},\quad e^{\pm i\pi x}\,\ln^3{x},\dots\}$. The set $\Phi_5=\Phi_1\cup \Phi_2\cup \Phi_3\cup
\Phi_4$ is also linearly
$v$-independent. 
\end{lemma}

\Proof We shall prove that $\Phi_2$ is linearly $v$-independent. For
the sake of convenience we denote $l_1(x)=\ln{x},\; l_2(x)=\ln(\ln{x}),\; l_3(x)=\ln(\ln(\ln{x})),
\dots$. So, we have $\Phi_2=\{l_n(x) : n=0, 1, 2,\dots\}$. It suffices to prove
that if $\sum_{k=1}^n\alpha_k\,l_k$ is a non-zero linear combination, then
$v\left(\sum_{k=1}^n\alpha_k\,l_k\right) = 0$. First note that this is
obviously true in the case $n=1$, since each $l_n(x)$ has valuation $0$. If
the sum has more than one term, we just need to show that $\sum_{k=1}^n\alpha_k\,l_k$ has a ``dominant
term,'' as $x\to\infty$. (Technically, this is the greatest term in the polynomial with respect
to the lexicographic ordering on the sequence of its exponents.) 
For this it suffices to show that  $m
> k$ implies $\lim_{x\to\infty}[(\lambda_m(x))^p/(\lambda_k(x))^q]=0$ whenever $p>0$ and $q>0$, which is
easily verified by l'Hopital's Rule. The rest of the sets are treated similarly  and we leave the
verification to the reader.
$\blacktriangle$ 

\begin{remark}[Linear Independence vs. Linear $v$-Independence] We observe that the sets in the above
four examples are also linearly independent over \C\, (in the usual sense). Notice that the set
$\Phi=\{\ln{x},\; 2\ln{x}\}$ is also   linearly $v$-independent but it is clearly linearly
dependent.  Conversely, the functions 
$f(x)=1$ and\, $g(x)=-1+1/x$ are linearly independent over \C, but $f$ and $g$ are not
$v$-independent, since $v(f+g) = v(1/x) = 1$ (not\, $=0$). We also observe that if a set $\Phi$ is
linearly $v$-independent over \C, then the set $\Phi\cup\{0\}$ is also linearly $v$-independent
over \C\, (in sharp contrast to the case of linear independence). 
\end{remark}
%_____________________________________________
\section{Pseudostandard Part}\label{S: Pseudostandard Part}

	We prove the existence of a particular type of linear homomorphism $\qst$ from $\mathcal{F}_v$ into
$\mathcal{F}_v$, with range in
$\mathcal{C}_v\cup\{0\}$, which is an extension of $\lim_{x\to\infty}$ (considered
also as a linear homomorphism). The applications of this construction appear in the next section.

\begin{definition}[Maximal Vector Spaces] \label{D: Maximal Vector Spaces}  Let $\mathcal{C}$ be a vector
subspace of $\mathcal{M}$. We say that
$\mathcal{C}$ is {\bf maximal in}
$\mathcal{C}_v$\; if\; $\C\subset\mathcal{C}\subset\mathcal{C}_v\cup\{0\}$ and there is 
no a vector subspace
$\mathcal{K}$ of $\mathcal{M}$ such that $\mathcal{C}\subsetneqq
\mathcal{K}\subset\mathcal{C}_v\cup\{0\}$.
\end{definition}
\begin{lemma}[Existence] \label{L: Existence}  Let $\Phi$ be
a linearly $v$-independent subset of
$\mathcal{M}$ (Definition~\ref{D: v-Independence}). Then there exists a vector subspace\, $\mathcal{C}$ of
$\mathcal{M}$ which is maximal in $\mathcal{C}_v$ and which contains $\Phi$.
\end{lemma}
\Proof Let $\mathcal{U}(\Phi)$ denote the set of all vector subspaces $\mathcal{K}$ of $\mathcal{M}$
such that $\mathbb{C}\cup\Phi\subseteq \mathcal{K}\subset\mathcal{C}_v\cup\{0\}$ and let
$\mathcal{U}(\Phi)$ be ordered by inclusion. Notice that $\mathcal{U}(\Phi)$ is non-empty since
$\Span(\Phi)\in\mathcal{U}(\Phi)$. We observe that every monotonic subset (chain)
$\mathcal{B}$ of
$\mathcal{U}(\Phi)$ is bounded from above by
$\bigcup_{B\in\mathcal{B}}\; B$. It follows that $\mathcal{U}(\Phi)$ has a maximal element
$\mathcal{C}$ (as desired), by Zorn's lemma. $\blacktriangle$\\

	We shall sometimes refer to $\Span(\Phi)$ as the ``explicit'' part of the space
$\mathcal{C}$ and to $\mathcal{C}\setminus\Span(\Phi)$ as the ``implicit'' part of
$\mathcal{C}$. 
\begin{theorem}[$v$-Completeness]\label{T: v-Completeness} Let $\mathcal{C}$ be a vector subspace
of
$\mathcal{M}$ which is maximal in
$\mathcal{C}_v$. Then: 

	{\bf (i)} $\mathcal{C}$ is {\bf $v$-complete} in $\mathcal{F}_v$ in the sense that: 
$\mu_v(\mathcal{C})=\mathcal{F}_v$, where
$\mu_v(\mathcal{C})$ is the $v$-monad of $\mathcal{C}$ ((\ref{E: Monad}), Section~\ref{S: Functions of
Moderate Growth}).

	{\bf (ii)} We have $\mathcal{F}_v=\mathcal{C}\oplus \mathcal{I}_v$ in the sense that every
$f\in\mathcal{F}_v$ has a unique asymptotic expansion
$f=\varphi+d\varphi$, where
$\varphi\in \mathcal{C}$ and $d\varphi\in \mathcal{I}_v$ ((\ref{E: v-infinitesimals}), Section~\ref{S:
Functions of Moderate Growth}).
\end{theorem}
\Proof (i) We have $\mu_v(\mathcal{C})\subseteq\mathcal{F}_v$, since $v(\varphi+d\varphi)\geq
\min\{v(\varphi), v(d\varphi)\}\geq 0$, by part (iii) of Theorem~\ref{T: Properties of v}. To show
$\mu_v(\mathcal{C})\supseteq\mathcal{F}_v$, suppose (on the contrary) that there exists
$g\in\mathcal{F}_v\setminus\mu_v(\mathcal{C})$, i.e. $v(g)\geq 0$ and $v(g-\varphi)\leq 0$ for all
$\varphi\in \mathcal{C}$. By letting $\varphi=0$, we conclude that $v(g)=0$, i.e. $g\in\mathcal{C}_v$. It
follows
$v(g-\varphi)=0$ for all $\varphi\in \mathcal{C}$ (by part (iii)
of Theorem~\ref{T: Properties of v}, again), i.e.
$g-\varphi\in\mathcal{C}_v$ for all $\varphi\in \mathcal{C}$. Thus we have 
$\mathcal{C}\subsetneqq
\Span(\mathcal{C}\cup\{g\})\subset\mathcal{C}_v\cup\{0\}$, contradicting the 
maximality of $\mathcal{C}$. 

	(ii) The existence of the asymptotic expansion $f=\varphi+d\varphi$ follows from (i). To
show the uniqueness, suppose that $\varphi+d\varphi=0$. It follows
$v(\varphi)=v(-d\varphi)=v(d\varphi)> 0$ implying $\varphi= d\varphi=0$ since $\mathcal{C}\subset
\mathcal{C}_v\cup\{0\}$. $\blacktriangle$\\

	The above result justifies the following definition.

\begin{definition}[Pseudostandard Part Mapping] \label{D: Pseudostandard Part Mapping} Let
$\mathcal{C}$ be a vector subspace of $\mathcal{M}$ which is maximal in
$\mathcal{C}_v$. We define $\qst_\mathcal{C}: \mathcal{F}_v\to \mathcal{F}_v$, by
$\qst_\mathcal{C}(\varphi+d\varphi)=\varphi$, where $\varphi\in \mathcal{C},\; d\varphi\in \mathcal{I}_v$.
We say that $\qst_\mathcal{C}$ is a {\bf pseudostandard part mapping} determined by $\mathcal{C}$. We shall
sometimes write simply $\qst$ instead of the more precise $\qst_\mathcal{C}$, suppressing the dependence
on $\mathcal{C}$, when the choice of $\mathcal{C}$ is clear from the context.
\end{definition}
\begin{theorem}[Properties of $\qst$]\label{T: Properties} Let $\mathcal{C}$ be a vector subspace of
$\mathcal{M}$ which is maximal in $\mathcal{C}_v$. Then:

	{\bf (i)} The pseudostandard part mapping $\qst_\mathcal{C}$ is a linear homomorphism from $\mathcal{F}_v$
into $\mathcal{F}_v$ with range $\ran{(\qst_\mathcal{C})}=\mathcal{C}$, i.e. $\qst_\mathcal{C}(\alpha
f+\beta g)=\alpha\qst_\mathcal{C}(f)+\beta\qst_\mathcal{C}(g)$, for all $f, g\in \mathcal{F}_v$ and all
$\alpha, \beta\in\mathbb{C}$.

	{\bf (ii)}  $\mathcal{C}$ consists of the fixed points of $\qst_\mathcal{C}$ in $\mathcal{F}_v$, i.e.
$\mathcal{C} = \{f\in\mathcal{F}_v \mid \qst_\mathcal{C}(\varphi)=\varphi \}$.
Consequently, we have $\qst_\mathcal{C}\circ\qst_\mathcal{C}= \qst_\mathcal{C}$. 

	{\bf (iii)} Suppose (in addition to the above) that $\Phi\subset\mathcal{C}$ for some linearly
$v$-independent subset $\Phi$ of $\mathcal{M}$ (Definition~\ref{D: v-Independence}).
Then we have $\qst_\mathcal{C}(\varphi+d\varphi)=\varphi$ for all $\varphi\in\Span(\Phi)$ and all
$d\varphi\in\mathcal{I}_v$. Consequently, the restriction $\qst_\mathcal{C}|{\C\oplus\mathcal{I}_v}$ is
a (proper) standard part mapping (in the sense of nonstandard analysis -- A. H. Lightstone and A.
Robinson~\cite{LiRob}), i.e. $\qst_\mathcal{C}(c+d\varphi)=c$ for all $c\in\C$ and all
$d\varphi\in\mathcal{I}_v$ (see the discussion after  Lemma~\ref{L: Connection with O's Symbols}). Or,
equivalently, $\qst_\mathcal{C}|{\C\oplus\mathcal{I}_v}$ coincides with $\lim_{x\to\infty}$ in the sense
that for every $f\in \C\oplus\mathcal{I}_v$ we have:
\[
\qst_\mathcal{C}(f)=\lim_{x\to\infty}f(x),
\]
\end{theorem}
\Proof (i) We have $f=\varphi+d\varphi$ and $g=\psi+d\psi$ for some $\varphi, \psi\in\mathcal{C},\;
d\varphi,\, d\psi\in\mathcal{I}_v$, by Theorem~\ref{T: v-Completeness}. Also
$\alpha f+\beta g=\alpha\varphi+\beta\psi+ \alpha d\varphi+ \beta d\psi$ and it is clear that
$\alpha\varphi+\beta\psi\in\mathcal{C}$ (since
$\mathcal{C}$ is a vector space) and
$\alpha d\varphi+ \beta d\psi\in\mathcal{I}_v$ (since $\mathcal{I}_v$ is an ideal in
$\mathcal{F}_v$, by Theorem~\ref{T: Rings and Ideals}). Thus
$\qst_\mathcal{C}(\alpha f+\beta
g)=\alpha\varphi+\beta\psi=\alpha\qst(\varphi)+\beta\qst(\psi)=
\alpha\qst_\mathcal{C}(f)+\beta\qst_\mathcal{C}(g)$,
as required, since  $\qst_\mathcal{C}(\varphi)=\varphi$ and  $\qst_\mathcal{C}(\psi)=\psi$, by the
definition of $\qst_\mathcal{C}$.

	(ii) As we mentioned already, $\qst_\mathcal{C}(\varphi)=\varphi$ for all $\varphi\in\mathcal{C}$, by the
definition of $\qst_\mathcal{C}$. Conversely, suppose that $\qst_\mathcal{C}(f)=f$ for some
$f\in\mathcal{F}_v$. It follows $f\in\mathcal{C}$ by the uniqueness of the asymptotic expansion
$f=\varphi+d\varphi$ (part (ii) of Theorem~\ref{T: v-Completeness}).

	(iii) follows from the definition of $\qst_\mathcal{C}$ taking into account that $\C\cup
\Phi\subset\mathcal{C}$.$\blacktriangle$
\begin{example}[The Case $\Phi=\Phi_5$]\label{Ex: The Case Phi=Phi5} Let $\mathcal{C}$ be
a vector subspace of $\mathcal{M}$ which is maximal in $\mathcal{C}_v$ and which contains
$\Phi_5$, where $\Phi_5$ is the set defined in Lemma~\ref{L: v-Independent Sets}
(Section~\ref{S: Linearly v-Independent Sets}). {\bf For the sake of simplicity we shall write
simply $\qst$ instead of the more precise $\qst_\mathcal{C}$}. By the above theorem, we have
$\qst(\varphi+d\varphi)=\varphi$ for all $\varphi\in\mathcal{C}$ (hence, for all
$\varphi\in\Span(\Phi_5)$) and all $d\varphi\in\mathcal{I}_v$. In particular, for every
$d\varphi\in\mathcal{I}_v$ and every $c\in\C$ we have:
$\qst(c+d\varphi)=c$; $\qst(c+1/\ln{x})=c+1/\ln{x}$,\; $\qst(\sin{x}+d\varphi)=\sin{x}$,\;
$\qst(\cos{x}+d\varphi)=\cos{x}$,\; 
$\qst(\sin{x}\cos{x}+d\varphi)=\sin{x}\cos{x}$,\; $\qst(\ln{x}+d\varphi)=\ln{x}$,\;
$\qst(\ln{(\ln{x})}+d\varphi)=\ln{(\ln{x})}$,\;
$\qst(\ln^n{x}+d\varphi)=\ln^n{x}$,\; $\qst(e^{\pm i\pi x}\ln^n{x}+d\varphi)=e^{\pm i\pi x}
\ln^n{x}$,\quad and $\qst\left((\sin{x})\, e^{\pm i\pi
x}\ln^n{x}+d\varphi\right)= (\sin{x})\, e^{\pm i\pi x} \ln^n{x},\quad n=0, 1, 2,
\dots$.
\end{example}
%-------------------------------------------------------------------------------------
\section{$v$-Asymptotic Series}\label{S: v-Asymptotic Series}

We consider a type of infinite series in $\mathcal{M}$ which we call a $v$-asymptotic series and
illustrate this concept by examples.

	Recall that a series $\sum_{n=0}^\infty f_n$ is called {\bf asymptotic} as $x\to\infty$, if 
$f_{n+1}=o(f_n)$ for all $n$. Notice that if an asymptotic series $\sum_{n=0}^\infty f_n$ is in
$\mathcal{M}$ (i.e. $f_n\in \mathcal{M}$ for all $n$), then we have $v(f_n) \leq v(f_{n+1})$ 
for all $n$. In contrast, we have the following similar concept:
\begin{definition}[$v$-Asymptotic Series]\label{D: v-Asymptotic Series} Let $F=\sum_{n=0}^\infty
f_n$ be a series in $\mathcal{M}$. Then:

	{\bf (i)} $F$ is called a {\bf $v$-asymptotic series}, as
$x\to\infty$, if it has the following two properties: 

	{\bf (a)} The set $\{x^{r_n}f_n\mid n\geq 0\}$ is
linearly $v$-independent over \C\, (Section~\ref{S: Linearly v-Independent Sets}), where
$r_n=v(f_n)$ and we let (by convention)
$x^\infty\cdot 0=0$.

	{\bf (b)} The restriction of the sequence $(r_n)$ on the set $\{n\mid f_n\not=0\}$ is a
strictly increasing (finite or infinite) sequence in \R. 

	Every $v$-asymptotic series $F$ can be written uniquely in its {\bf
canonical form}:
\begin{equation}\label{E: Canonical}
F=\sum_{n=0}^\infty\; \frac{\varphi_n(x)}{x^{r_n}},
\end{equation}
where $\varphi_n(x)=x^{r_n}f_n(x)$ and we let (by
convention) $\frac{0}{x^\infty}=0$. We call the sets $Supp(F)=\{n\mid f_n\not=0\}=
\{n\mid \varphi_n\not=0\}$,\; $Exp(F)=\{r_n\mid n\geq 0\}$ and $Coe\!f(F)=\{x^{r_n}f_n\mid n\geq
0\}=\{\varphi_n\mid n\geq 0\}$ the {\bf support}, the {\bf set of exponents} and the {\bf set of
coefficients} of $F$, respectively. 
\end{definition}

	In what follows we shall often {\bf assume that the $v$-asymptotic series are already presented in their
canonical form} (\ref{E: Canonical}). Notice that
for the $v$-asymptotic series we have $Coe\!f(F)\subset \mathcal{C}_v\cup\{0\}$ 
(Section~\ref{S: Linearly v-Independent Sets}).
\begin{notation}\label{N: Sets of Series} Let $K$ be a subset of $\mathcal{C}_v\cup\{0\}$.

	{\bf (i)} We denote by $K(1/x^\omega)$ the set of
all $v$-asymptotic series in the form (\ref{E: Canonical}) with coefficients $\varphi_n$ in $K$.

	{\bf (ii)} We denote by $K\bra 1/x\ket$ the set of
all $v$-asymptotic series in the form (\ref{E: Canonical}) with coefficients $\varphi_n$ in $K$ and
such that $\lim_{n\to\infty}r_n=\infty$.
\end{notation}

	We should notice that $K(1/x^\omega)$ consists of all Hahn's series \cite{hH} with
coefficients in $K$ and with a set of exponents which is an increasing sequence is $\R$ (hence, the origin
of
$\omega$ in the notation). We
observe as well that in the particular case
$K=\mathcal{C}_v\cup\{0\}$, the set $K(1/x^\omega)$ coincides with the set of all $v$-asymptotic series
in $\mathcal{M}$, while $K\bra 1/x\ket$ coincides with the set of all $v$-asymptotic series
(\ref{E: Canonical}) in
$\mathcal{M}$ with $\lim_{n\to\infty}r_n=\infty$.

	Here are several examples:

\begin{example} Let $c_n$ be a (arbitrary) sequence in \C\, and $m\in\Z$. Then the series
$\sum_{n=0}^{\infty}\; \frac{c_n}{x^{m+n}}$ is both asymptotic and $v$-asymptotic. 
\end{example}
\begin{example}The series: $\sum_{n=0}^\infty \frac{\ln^n{x}}{x^n}$\; $\sum_{n=0}^\infty \frac{e^{i n
x}}{x^n}$,\; $\sum_{n=1}^\infty\sqrt[n]{x}$,\; $\sum_{n=1}^\infty\; \ln^n{x}$ and $\ln{x}+\ln(\ln{x}) +
\ln(\ln(\ln{x})))+\dots$ are all asymptotic. All but the last two are also
$v$-asymptotic. The last two series are not
$v$-asymptotic, since 
$v(\ln^n{x})=0$ and $v(\ln{x})= v(\ln(\ln{x})) = v(\ln(\ln(\ln{x}))))=\dots = 0$. 
\end{example}

\begin{example}\label{Ex: Sine Series} In contrast to the above the Sine-series: 
\begin{align}\notag
\frac{\pi}{2} & +
\cos{x}\left(-\frac{1}{x}+\frac{2!}{x^3}-\frac{4!}{x^5}+\dots+\frac{(2n)!\,(-1)^{n+1}}{x^{2n+1}}+\dots\right)+\\
&+\sin{x}\left(-\frac{1}{x^2}+\frac{3!}{x^4}-\frac{5!}{x^6}+\dots+\frac{(2n+1)!(-1)^{n+1}}{x^{2n+2}}+\dots\right),\notag
\end{align}
is $v$-asymptotic (see the set $\Phi_1$ in Lemma~\ref{L: v-Independent Sets}, Section~\ref{S:
Linearly v-Independent Sets}), but not asymptotic since

\[
\frac{(\sin{x})(2n+1)!(-1)^{n+1}}{x^{2n+2}}\neq o\left(\frac{(\cos{x})(2n)!\,(-1)^{n+1}}{x^{2n+1}}\right).
\]
\end{example}
\begin{example} The series $\sum_{n=0}^\infty\frac{1+n+1/\sqrt[n]{x}}{x^n}$ is asymptotic,
but not $v$-asymptotic since the set of the coefficients $\{1+n+\frac{1}{\sqrt[n]{x}}\mid n\geq 0\}$ is
not linearly $v$-independent over $\mathbb{C}$.

\end{example}

\begin{example} Finally, the series $\sum_{n=0}^\infty \frac{x^n}{n!}$ and $\sum_{n=0}^\infty
\frac{1}{\sqrt[n]{x}}$ are neither asymptotic, nor $v$-asymptotic as $x\to\infty$ 
(the first is convergent in \R,  while the second is divergent in \R,
but the convergence or divergence in \R\, is irrelevant to their asymptotic properties).
\end{example}
%_____________________________________
\section{$v$-Asymptotic Expansion of a Function}\label{S: v-Asymptotic Expansion of a Function}

	In this section we discuss the concept of the asymptotic expansion of a function of moderate growth at
infinity in a $v$-asymptotic series. 

\begin{definition} Let $f\in \mathcal{M}$ (Definition~\ref{D: Functions of Moderate Growth}) and let
$\sum_{n=0}^\infty\,\frac{\varphi_n(x)}{x^{r_n}}$ be a $v$-asymptotic series in $\mathcal{M}$
(Section~\ref{S: v-Asymptotic Series}). We say that the series
$\sum_{n=0}^\infty\,\frac{\varphi_n(x)}{x^{r_n}}$ is a {\bf $v$-asymptotic expansion} of $f(x)$ (or, that
$f(x)$ is an {\bf asymptotic sum} of\,  $\sum_{n=0}^\infty\,\frac{\varphi_n(x)}{x^{r_n}}$), as
$x\to\infty$, in
symbol, $f(x)\leadsto
\sum_{n=0}^\infty\,\frac{\varphi_n(x)}{x^{r_n}}$, if for every integer 
$n\geq 0$ for which $\varphi_n\not=0$, we
have:
\begin{equation}\label{E: Residual}
 f(x) - \sum_{k=0}^n\, \frac{\varphi_k(x)}{x^{r_k}} = o(1/x^{r_n}),\quad \text{as}\; x\to\infty.
\end{equation}
\end{definition}
\begin{remark} Notice that $r_n=\infty$ \ifff $\varphi_n=\varphi_{n+1}=\varphi_{n+2}=\dots=0$ 
\ifff
$f(x)-\sum_{k=0}^n\,\frac{\varphi_k(x)}{x^{r_k}}=f(x)-\sum_{k=0}^{n-1}\,
\frac{\varphi_k(x)}{x^{r_k}}\in\mathcal{N}$ since, by convention, $\frac{0}{x^\infty}=0$.
\end{remark}
\begin{lemma}\label{L: Characterization of AsyExp} Let $f\in
\mathcal{M}$ and 
$\sum_{n=0}^\infty\,\frac{\varphi_n(x)}{x^{r_n}}$ be $v$-asymptotic series in
$\mathcal{M}$. Then the following are
equivalent:

	{\bf (i)}\quad $f(x) \leadsto \sum_{n=0}^\infty\,\frac{\varphi_n(x)}{x^{r_n}}$.

	{\bf (ii)}\, $x^{r_n}\left[f(x) -\sum_{k=0}^n\,\frac{\varphi_k(x)}{x^{r_k}}\right]\in\mathcal{I}_v$ for
all $n\geq 0,\, \varphi_n\not=0$.

	{\bf (iii)} $x^{r_n}\left[f(x)
-\sum_{k=0}^{n-1}\,\frac{\varphi_k(x)}{x^{r_k}}\right]\in\mathcal{C}_v$ for all $n\geq
0,\, \varphi_n\not=0$.
\end{lemma}
\Proof (i) $\Rightarrow$ (ii): Suppose, first, that $\varphi_{n+1}\not=0$. We have $f(x) -
\sum_{k=0}^n\,
\frac{\varphi_k(x)}{x^{r_k}} =
\frac{\varphi_{n+1}(x)}{x^{r_{n+1}}}+o(1/x^{r_{n+1}})$. It follows 
\begin{align}\notag
&v\left(x^{r_n}\left[f(x)
-\sum_{k=0}^n\,\frac{\varphi_k(x)}{x^{r_k}}\right]\right)=v\left(\frac{\varphi_{n+1}(x)+o(1)}
{x^{r_{n+1}-r_n}}\right)=\\\notag
&=v\left(\varphi_{n+1}(x)+o(1)\right)+r_{n+1}-r_n>0,
\end{align}
as required, by Theorem~\ref{T: Properties of v}, since $v\left(\varphi_{n+1}(x)+o(1)\right)\geq 0$ and
$r_{n+1}-r_n>0$ (by assumption). If $\varphi_{n+1}=0$, then
$f(x)-\sum_{k=0}^n\,\frac{\varphi_k(x)}{x^{r_k}}\in\mathcal{N}$ implying $x^{r_n}\left[f(x)
-\sum_{k=0}^n\,\frac{\varphi_k(x)}{x^{r_k}}\right]\in\mathcal{N}$ and (ii) follows since
$\mathcal{N}\subset\mathcal{I}_v$.

	(i) $\Leftarrow$ (ii) follows directly by part (i) of Lemma~\ref{L: Connection
with O's Symbols}.

	(ii) $\Rightarrow$ (iii) We denote $x^{r_{n}}\left[f(x)
-\sum_{k=0}^n\,\frac{\varphi_k(x)}{x^{r_k}}\right]=d\varphi_n(x)$ and observe that\newline
$v\left(x^{r_n}\left[f(x)
-\sum_{k=0}^{n-1}\,\frac{\varphi_k(x)}{x^{r_k}}\right]\right)=v\left(\varphi_n(x)+d\varphi_n(x)\right)=0$,
as required, since $d\varphi_n\in\mathcal{I}_v$, by assumption. 

	(ii) $\Leftarrow$ (iii): Suppose, first, that $\varphi_{n+1}\not=0$. We denote $x^{r_{n+1}}\left[f(x)
-\sum_{k=0}^n\,\frac{\varphi_k(x)}{x^{r_k}}\right]=\psi_n(x)$ and observe that $\psi_n\in
\mathcal{C}_v$, by assumption. It follows
\[
x^{r_n}\left[f(x)
-\sum_{k=0}^n\,\frac{\varphi_k(x)}{x^{r_k}}\right]=\frac{\psi_n(x)}{x^{r_{n+1}-r_n}}\in\mathcal{I}_v.
\]
since $r_{n+1}-r_n>0$. If $\varphi_{n+1}=0$, then
$f(x)-\sum_{k=0}^n\,\frac{\varphi_k(x)}{x^{r_k}}\in\mathcal{N}$ implying \newline 
$x^{r_n}\left[f(x)
-\sum_{k=0}^n\,\frac{\varphi_k(x)}{x^{r_k}}\right]\in\mathcal{N}$ and (ii) follows.
$\blacktriangle$
\begin{theorem}[Convergent Series in $\mathcal{M}$]\label{T: Convergent Series in M} Let
$f\in\mathcal{M}$ and $\sum_{n=0}^\infty\,\frac{\varphi_n(x)}{x^{r_n}}$ be a $v$-asymptotic series
in $\mathcal{M}$ such that $\lim_{n\to\infty}r_n=\infty$. Then the following are equivalent:
	
	{\bf (i)} $f(x)\leadsto\sum_{n=0}^\infty\,\frac{\varphi_n(x)}{x^{r_n}}$.

	{\bf (ii)} $f(x)=\sum_{n=0}^\infty\,\frac{\varphi_n(x)}{x^{r_n}}$ in the sense that
$f(x)=\lim_{n\to\infty}\; \sum_{k=0}^n\, \frac{\varphi_k(x)}{x^{r_k}}$ in $(\mathcal{M},
d_v)$ (Definition~\ref{D: Pseudometric}, Section~\ref{S: Functions of Moderate Growth}).
\end{theorem}
\Proof $(i)\Rightarrow (ii)$: We have $f(x)-\sum_{k=0}^n\, \frac{\varphi_k(x)}{x^{r_k}}=o(1/x^{r_n})$
for all
$n$ which implies (ii), by Corollary~\ref{C: For LC} (Section~\ref{S: Functions of Moderate Growth}),
since
$\lim_{n\to\infty}r_n=\infty$, by assumption. 

	$(i)\Leftarrow(ii)$: Suppose, first, that $f\in\mathcal{N}$ (in particular, $f=0$). It follows \newline
$\lim_{n\to\infty} v\left(\sum_{k=0}^n\,
\frac{\varphi_k(x)}{x^{r_k}}\right)=\infty$, by part~(i) of Theorem~\ref{T: Convergence in M}. It follows
$r_0=r_1=r_2=\dots=\infty$ and $\varphi_0(x)=\varphi_1(x)=\varphi_2(x)=\dots=0$, thus, (i) hods (in the
form $f(x)\leadsto 0+0+\dots$). Suppose that $f\notin\mathcal{N}$. It follows $v(f)=v\left(\sum_{k=0}^n\,
\frac{\varphi_k(x)}{x^{r_k}}\right)$ for all sufficiently large $n$, by part~(ii) of Theorem~\ref{T:
Convergence in M} (Section~\ref{S: Functions of Moderate Growth}). Thus
$v(f)=r_0$ and more generally, $v\left(f(x)-\sum_{k=0}^n\,\frac{\varphi_k(x)}{x^{r_k}}\right)=r_{n+1}$
for all $n$. Similarly, we observe that $\qst_\mathcal{C}\left(x^{r_0}f(x)\right)=\varphi_0(x)$
and
$\qst_\mathcal{C}\left(x^{r_{n+1}}\,\left[f(x)-\sum_{k=0}^{n}\frac{\varphi_k(x)}{x^{r_k}}\right]\right)=\varphi_{n+1}(x)$
for all $n$, where
$\mathcal{C}$ is a maximal vector space containing the set $\{\varphi_n \mid  n\in\N\}$. It follows that
for each
$n$ we have $f(x)-\sum_{k=0}^{n}\frac{\varphi_k(x)}{x^{r_k}}\leadsto
\sum_{k=n+1}^{\infty}\frac{\varphi_k(x)}{x^{r_k}}$ which implies (i), as required. 
$\blacktriangle$
%________________________________________
\section{Existence and Uniqueness Result}\label{S: Existence and Uniqueness Result}

The purpose of this section is to show that every function with moderate growth at infinity 
has a $v$-asymptotic expansion in $\mathcal{M}$. 
\begin{theorem}[Existence and Uniqueness]\label{T: Existence and Uniqueness} Every function with
moderate growth at infinity  has a $v$-asymptotic expansion in $\mathcal{M}$. More precisely, 
let $\mathcal{C}$ be a vector subspace of $\mathcal{M}$ which is maximal in $\mathcal{C}_v$
(Definition~\ref{D: Maximal Vector Spaces}, Section~\ref{S: Pseudostandard Part}). Then every
function
$f$ in
$\mathcal{M}$ has a unique
$v$-asymptotic expansion (as $x\to\infty$):
\begin{equation}\label{E: AsyExp}
f(x) \leadsto \sum_{n=0}^\infty\,\frac{\varphi_n(x)}{x^{r_n}},
\end{equation}
with coefficients $\varphi_n$ in $\mathcal{C}$ (that is, within the set $\mathcal{C}(1/x^\omega)$,
Notation~\ref{N: Sets of Series}, Section~\ref{S: v-Asymptotic Series}). 
\end{theorem}
\Proof	{\em (A) Existence:} Our first goal is to show the existence of a
$v$-asymptotic series in (\ref{E: AsyExp}). Let
$\qst_\mathcal{C}:
\mathcal{F}_v\to\mathcal{F}_v$ be the
pseudostandard part mapping determined by $\mathcal{C}$ (Definition~\ref{D: Pseudostandard Part Mapping},
Section~\ref{S: Pseudostandard Part}). The function $f(x)$ determines the
sequences $(\varphi_n)$ and $(r_n)$ by the following {\bf recursive formulas:}
\begin{align}\label{E: Recursive}
&r_0=v(f),  &&\varphi_0(x)=\qst_\mathcal{C}\left(x^{r_0}f(x)\right),\\ \notag
&r_1=v\left(f(x)-\frac{\varphi_0(x)}{x^{r_0}}\right),&&\varphi_1(x)=\qst_\mathcal{C}
\left(x^{r_1}\left[f(x)-\frac{\varphi_0(x)}{x^{r_0}}\right]\right),\\\notag
&\dots,\text{etc.},\\
&r_n=v\left(f(x)-\sum_{k=0}^{n-1}\frac{\varphi_k(x)}{x^{r_k}}\right), &&
\varphi_n(x)=\qst_\mathcal{C}\left(x^{r_n}
\left[f(x)-\sum_{k=0}^{n-1}\frac{\varphi_k(x)}{x^{r_k}}\right]\right),\label{E: n-Recursive}
\end{align}
\begin{align}
&r_{n+1}=v\left(f(x)-\sum_{k=0}^{n}\frac{\varphi_k(x)}{x^{r_k}}\right),
&&\varphi_{n+1}(x)=\qst_\mathcal{C}\left(x^{r_{n+1}}\,
\left[f(x)-\sum_{k=0}^{n}\frac{\varphi_k(x)}{x^{r_k}}\right]\right),\notag
\end{align}
\noindent $n=0, 1, 2,\dots$, where, by convention, we let $x^\infty\cdot
\varphi=0$ whenever $\varphi\in\mathcal{N}$. Next, we observe that $\varphi_n\in\mathcal{C}$ (since
$\ran{(\qst_\mathcal{C})}=\mathcal{C}$) and that $(r_n)$ is a non-decreasing sequence in
$\R\cup\{\infty\}$. Our next goal is to show that the series on the RHS of (\ref{E:
AsyExp}) is a $v$-asymptotic series. Suppose, first, that $r_n=\infty$ for some $n$. It follows
$f(x)-\sum_{k=0}^{n-1}\frac{\varphi_k(x)}{x^{r_k}}\in\mathcal{N}$ thus
$x^\infty\cdot[f(x)-\sum_{k=0}^{n-1}\frac{\varphi_k(x)}{x^{r_k}}]=0$ (by convention), implying
that
$\varphi_n=0$ (since $\mathcal{N}\subset\mathcal{I}_v$).
Consequently, it follows that
$r_k=\infty$ and
$\varphi_k=0$ for all $k\geq n$, as required (which means that in this case the series is a finite sum).
Suppose now that $r_{n+1}\not=\infty$ for some $n$. It follows that $r_n\not=\infty$ as well, by what was
just proved above. We have to show that $r_n<r_{n+1}$. We observe that the second formula in (\ref{E:
n-Recursive}) implies
\begin{equation}\label{E: n-th}
f(x)=\sum_{k=0}^{n}\frac{\varphi_k(x)}{x^{r_k}}+\frac{d\varphi_n(x)}{x^{r_n}}.
\end{equation}
for some
$d\varphi_n(x)\in\mathcal{I}_v$. It follows
\[\notag
r_{n+1}=v\left(f(x)-\sum_{k=0}^{n}\frac{\varphi_k(x)}{x^{r_k}}\right)=
v\left(\frac{d\varphi_n(x)}{x^{r_n}}\right)\geq
v(d\varphi_n)+ v(1/x^{r_n})>r_n,
\]
by Theorem~\ref{T: Properties of v}, since $v(d\varphi_n)>0$ and $v(1/x^{r_n})=r_n$. Next, we observe
that the set $\Phi=\{\varphi_n \mid n\geq 0\}$ is linearly $v$-independent over
\C, by Lemma~\ref{L: Trivial} (Section~\ref{S: Linearly
v-Independent Sets}). Thus the series on the RHS of (\ref{E:
AsyExp}) is a
$v$-asymptotic series. To show that the asymptotic expansion (\ref{E: AsyExp}) holds, suppose that
$\varphi_n\not=0$ (hence, $r_{n}\not=\infty$). The formula (\ref{E: n-th}) implies $x^{r_n}
\left[f(x)-\sum_{k=0}^{n}\frac{\varphi_k(x)}{x^{r_k}}\right]=d\varphi_n(x)$ which implies
(\ref{E: AsyExp}), by Lemma~\ref{L: Characterization of AsyExp}, since
$d\varphi_n(x)\in \mathcal{I}_v$. This proves the existence of a $v$-asymptotic expansion of $f$
within the set of series $\mathcal{C}(1/x^\omega)$.

	{\em (B) Uniqueness:} Suppose that 
\begin{equation}
f(x) \leadsto \sum_{n=0}^\infty\,\frac{\psi_n(x)}{x^{s_n}},\quad \text{as\;} x\to\infty,
\end{equation} 
for another $v$-asymptotic series
$\sum_{n=0}^\infty\,\frac{\psi_n(x)}{x^{s_n}}$, where $\psi_n\in\mathcal{C}$. Suppose as well that
$\psi_n\not=0$  (hence, $s_{n}\not=\infty$) for some $n$. It follows that
$x^{s_n}\left[f(x) -\sum_{k=0}^{n-1}\,\frac{\psi_k(x)}{x^{s_k}}\right]\in\mathcal{C}_v$, by part~(iii) of
Lemma~\ref{L: Characterization of AsyExp}. Thus 
\begin{equation}\label{E: Valuation}
s_n= v\left(f(x)-\sum_{k=0}^{n-1}\,\frac{\psi_k(x)}{x^{s_k}}\right),
\end{equation}
by part~(ii) of Theorem~\ref{T: Properties of v}. Similarly, we have $x^{s_n}\left[f(x)
-\sum_{k=0}^{n}\,\frac{\psi_k(x)}{x^{s_k}}\right]\in\mathcal{I}_v$, by part~(ii) of
Lemma~\ref{L: Characterization of AsyExp}. It follows 
\begin{equation}\label{E: Coefficients}
\psi_n(x)=
\qst_\mathcal{C}\left(x^{s_n}\left[f(x)-\sum_{k=0}^{n-1}\,\frac{\psi_k(x)}{x^{s_k}}\right]\right),
\end{equation}
since $\qst_\mathcal{C}(\psi_n)=\psi_n$. Since $\psi_k\not=0$ for all $k\leq n$ we can start from
$k=0$. Formula (\ref{E: Valuation}) gives
$s_0=v(f(x))=r_0$ and formula (\ref{E: Coefficients}) gives $\psi_0(x)=
\qst_\mathcal{C}\left(x^{s_0}f(x)\right)=\qst_\mathcal{C}\left(x^{r_0}f(x)\right)=\varphi_0(x)$.
Similarly, $k=1$ gives
\begin{align}\notag
&s_1=v\left(f(x)-\frac{\psi_0(x)}{x^{s_0}}\right)=v\left(f(x)-\frac{\varphi_0(x)}{x^{r_0}}\right)=r_1,\\\notag
&\psi_1(x)=\qst_\mathcal{C}\left(x^{s_1}\left[f(x)-\frac{\psi_0(x)}{x^{s_0}}\right]\right)=\qst_\mathcal{C}\left(x^{r_1}\left[f(x)-
\frac{\varphi_0(x)}{x^{r_0}}\right]\right)=\varphi_1(x).
\end{align}
Continuing we obtain $s_k=r_k$ and $\psi_k(x)=\varphi_k(x)$ for all $k\leq n$. Suppose, now, that
$\psi_n=0$ for some $n$ which implies $s_{n+1}=s_{n+2}=\dots=\infty$ and
$\psi_n=\psi_{n+1}=\psi_{n+2}=\dots=0$. Without loss of generality we can assume that
$\psi_{n-1}\not=0$, hence, $s_{n-1}\not=\infty$. So, we have $s_k=r_k$ and $\psi_k(x)=\varphi_k(x)$ for
all
$k\leq n-1$. It follows
\[
f(x)-\sum_{k=0}^{n-1}\,\frac{\varphi_k(x)}{x^{r_k}}=
             f(x)-\sum_{k=0}^{n-1}\,\frac{\psi_k(x)}{x^{s_k}}\in\mathcal{N},
\]
which, on its turn, implies $r_n=r_{n+1}=r_{n+2}=\dots=\infty$ and
$\varphi_n=\varphi_{n+1}=\varphi_{n+2}=\dots=0$. Summarizing, we have $s_k=r_k$ and $\psi_k=\varphi_k$
for all $k\geq 0$ which proves the uniqueness of the $v$-asymptotic expansion within
$\mathcal{C}(1/x^\omega)$ (Notation~\ref{N: Sets of Series}, Section~\ref{S: v-Asymptotic Series}). 
To complete the proof we observe that there exist
maximal linear spaces $\mathcal{C}$ in $\mathcal{M}$, by 
Lemma~\ref{L: Existence}, Section~\ref{S: Pseudostandard Part} (in particular, there exists
$\mathcal{C}$ such that $\Phi_5\subset\mathcal{C}$, where $\Phi_5$ is defined in Lemma~\ref{L:
v-Independent Sets}).
$\blacktriangle$

	The next corollary is written for those readers who do not feel comfortable with the concept of 
``maximal linear space''.
\begin{corollary}[On Uniqueness Again]\label{C: On Uniqueness Again} Let (\ref{E: AsyExp}) holds for
some moderate function 
$f\in\mathcal{M}$ and some $v$-asymptotic series $\sum_{n=0}^\infty\,\frac{\varphi_n(x)}{x^{r_n}}$.
Let $\mathcal{K}$ be a linear subspace of $\mathcal{M}$ such that
$\Phi\subseteq\mathcal{K}\subset\mathcal{C}_v\cup\{0\}$, where $\Phi=\{\varphi_n\mid n=0,
1,2,\dots\}$. Then (\ref{E: AsyExp}) is the only
$v$-asymptotic expansion of $f$ within the set of $v$-asymptotic series
$\mathcal{\mathcal{K}}(1/x^\omega)$ (Notation~\ref{N: Sets of Series}, Section~\ref{S: v-Asymptotic
Series}). 
\end{corollary}
\Proof Let $\mathcal{C}$ be a vector subspace of $\mathcal{M}$ which is maximal 
in $\mathcal{C}_v$ and which contains $\mathcal{K}$. 
Notice that the existence of $\mathcal{C}$ is guaranteed by 
Lemma~\ref{L: Existence} (Section~\ref{S: Pseudostandard Part}). The asymptotic expansion 
(\ref{E: AsyExp}) is unique in $\mathcal{\mathcal{C}}(1/x^\omega)$, by the above theorem. The latter
implies the uniqueness of (\ref{E: AsyExp}) in $\mathcal{\mathcal{K}}(1/x^\omega)$.
$\blacktriangle$

	In the next example we explain how to use Theorem~\ref{T: Existence and Uniqueness} in asymptotic
analysis.
%______________________________________
\begin{example}[An Integral with a Large Parameter]\label{Ex: An Example: An Integral 
with Large Parameter} Let $f\in C^\infty[-\pi, \pi]$ and define:
\begin{equation}\label{E: Integral}
I(x) = \int_{-\pi}^\pi e^{ixy-y^2\ln{x}} f(y)dy.    
\end{equation}
Integrals of the form (\ref{E: Integral}) arise in quantum
statistical mechanics and random matrix theory (M. L. Mehta~\cite{mlMehta}), where they are related to
the probability that no eigenvalue of a large  Hermitian matrix lies in the interval $(-x, x)$
(see also E.~Basor~\cite{eBasor99} and C.~P.~Hughes, J.~P.~Keating,
N.~O'Connell~\cite{HughesKeatingO'Connell01}). Since $I(x)$ is a function with moderate
growth at infinity (Section~\ref{S: Functions of Moderate Growth}), it follows that $I(x)$ has a
$v$-asymptotic expansion in $\mathcal{M}$, by Theorem~\ref{T: Existence and Uniqueness}.
Notice that Theorem~\ref{T: Existence and Uniqueness} is used to guarantee the existence 
of a $v$-asymptotic expansion in advance; before plunging into possibly hard and time consuming
calculations. Next, we note that the
recursive formulas (\ref{E: Recursive}) (used in the proof of Theorem~\ref{T: Existence and Uniqueness})
are rarely efficient for explicit asymptotic calculations (just as the recursive
formulas (\ref{E: Limits}) in the Introduction
are rarely efficient for explocit calculations). The
explicit calculations should be done by any of the methods known  to asymptotic analysis. In the case
of $I(x)$ the ``integrating by parts'' produces the following
$v$-asymptotic expansion:
\begin{equation}\label{E: IAsyExp}
I(x) \leadsto \sum_{n=0}^\infty\,\frac{\varphi_n(x)}{x^{\pi^{2}+1+n}},
\end{equation}
where the sequence $(\varphi_n)$ is determined by the recursive formulas:
\begin{align}\notag
&\varphi_0(x) = \frac{1}{i}\left[f(\pi)e^{i\pi x} - f(-\pi)e^{-i\pi x}\right],\\\notag
&\varphi_1(x) = f'(\pi)e^{i\pi x} - f'(-\pi)e^{-i\pi x}-(2\pi\ln{x})\left[f(\pi)e^{i\pi x} +
f(-\pi)e^{-i\pi x}\right], \notag
\end{align}
and, more generally,
\begin{align} \notag 
\varphi_n(x) &=\frac{(-1)^n}{i^{n+1}}\sum_{k=0}^{[\frac{n}{2}]}\binom{n}{2k}(2\pi\ln{x})^{2k}\left[f^{(n-2k)}(\pi)e^{i\pi x}-f^{(n-2k)}(-\pi)e^{-i\pi
x}\right]-\label{E: The coefficients} \\
 &- \frac{(-1)^n}{i^{n+1}}\sum_{k=0}^{[\frac{n-1}{2}]} \binom{n}{2k+1}(2\pi\ln{x})^{2k+1} \left[f^{(n-2k-1)}(\pi) e^{i\pi x} +
f^{(n-2k-1)}(-\pi)e^{-i\pi x}\right].\notag
\end{align}
Now is the time to establish our uniqueness result. First, we observe that
$\Span(\{\varphi_n\mid n=0, 1,2,\dots\})=\Span(\Phi_4)$, where $\Phi_4$ is the set of functions
defined in Lemma~\ref{L: v-Independent Sets} (Section~\ref{S: Linearly v-Independent Sets}). We
conclude that the series in (\ref{E: IAsyExp}) is a $v$-asymptotic series since $\Phi_4$ is a linearly
$v$-independent set and the sequence of exponents $r_n=\pi^{2}+1+n$ is strictly increasing.  Let
$\mathcal{K}$ be a vector subspace of
$\mathcal{M}$ such that 
$\Phi_4\subset\mathcal{K}\subset\mathcal{C}_v\cup\{0\}$
and let $\mathcal{K}(1/x^\omega)$ denote the set of all $v$-asymptotic series with
coefficients in $\mathcal{K}$ (Notation~\ref{N: Sets of Series}, Section~\ref{S: v-Asymptotic
Series}). We conclude that the
$v$-asymptotic expansion (\ref{E: IAsyExp}) of $I(x)$ is unique within the set of series
$\mathcal{K}(1/x^\omega)$, by Corollary~\ref{C: On Uniqueness Again}. In particular, let
$\mathcal{C}$ be a vector subspace of $\mathcal{M}$ which is maximal in $\mathcal{C}_v$ and which
contains $\Phi_4$. Then the
$v$-asymptotic expansion (\ref{E: IAsyExp}) of $I(x)$ is unique
within set of series $\mathcal{C}(1/x^\omega)$. Finally, it is curious to
observe that $I(x) =\sum_{n=0}^\infty\,\frac{\varphi_n(x)}{x^{\pi^{2}+1+n}}$, by Theorem~\ref{T:
Convergent Series in M}, since $\lim_{n\to\infty}(\pi^2+1+n)=\infty$. That means that the series
(\ref{E: IAsyExp}) is convergent and  $I(x)
=\lim_{n\to\infty}\,\sum_{k=0}^n\,\frac{\varphi_k(x)}{x^{\pi^{2}+1+k}}$ in the metric space
$(\mathcal{M}, d_v)$ (Definition~\ref{D: Pseudometric}, Section~\ref{S: Functions of Moderate
Growth}).
\end{example}

%____________________________________________
\section{Colombeau's Generalized Numbers}\label{S: Colombeau's Generalized Numbers}

	As an application to the nonlinear theory of generalized functions (J.F.
Colombeau~\cite{jfCol84}-\cite{jfCol85}), we show that every Colombeau's generalized number has a
$v$-asymptotic expansion. Since the Colombeau generalized functions can be characterized as pointwise
functions in the ring of Colombeau' generalized numbers (M. Kunzinger and M.
Oberguggenberger~\cite{KunzOber99}), it follows that every Colombeau's generalized function
(including every Schwartz distribution) has (pointwisely) a $v$-asymptotic expansion. This result has
an importance for finding weak generalized solutions of some partial differential equations with
singularities (M. Oberguggenberger~\cite{mOber92}).
\begin{definition}[Colombeau's Ring]\label{D: Colombeau's Ring}

	{\bf (i)} We call the factor-ring
$\overline{\C}=\mathcal{M}/\mathcal{N}$ Colombeau's ring of {\bf complex generalized numbers}
(Section~\ref{S: Functions of Moderate Growth}). We denote by
$q:\mathcal{M}\to\overline{\C}$ the corresponding quotient mapping. The Colombeau's ring of {\bf real
generalized numbers} is 
\[
\overline{\R}=\{q(f)\mid f\in\mathcal{M},\; f(x)\in\R \text{\; for all sufficiently
large\;} x\}.
\]
The real generalized number $\lambda=q(id)$ is called the {\bf scale} of $\overline{\C}$, where $id(x)=x$
for all $x\in\mathbb{R}$.

	{\bf (ii)} We define a pseudovaluation $v:
\overline{\C}\to\R\cup\{\infty\}$ and a metric\newline
$d_v:
\overline{\C}\times\overline{\C}\to\R$ inherited from $\mathcal{M}$, i.e. by\,
$v(q(f))=v(v)$ and
$d_v(q(f), v(g))=d_v(f, g)$, respectively. We denote by $(\overline{\C}, d_v)$ the corresponding
metric space.

	{\bf (iii)} We denote by $\mathcal{I}_v(\overline{\C}), \mathcal{C}_v(\overline{\C})$ and
$\mathcal{F}_v(\overline{\C})$ the subsets of $\overline{\C}$ with positive, zero and non-negative
pseudovaluation, respectively.

	{\bf (iv)} Let $C$ be a vector
subspace of $\overline{\C}$. We say that
$C$ is {\bf maximal in}
$\mathcal{C}_v(\overline{\C})$\; if\; $\mathbb{C}\subset C\subset\mathcal{C}_v(\overline{\C})\cup\{0\}$
and there is no a vector subspace
$\mathbb{K}$ of $\overline{\C}$ such that $C\subsetneqq
\mathbb{K}\subset\mathcal{C}_v(\overline{\C})\cup\{0\}$.

	{\bf (v)} We say that a set $\mathbb{F}\subset\overline{\C}$ is {\bf linearly $v$-independent over} \C\,
if
$v(a)=0$ for all
$a\in\Span(\mathbb{F}),\; a\not=0$. 

\end{definition}
\begin{remark} The ring $\overline{\C}$ defined above is, actually, isomorphic to
the original Colombeau ring of generalized numbers (under the mapping $x\to 1/\epsilon$), introduced in
(J.F. Colombeau~\cite{jfCol85}, \S 2.1). We shall ignore the difference between these two isomorphic
rings. Notice that both $\overline{\R}$ and
$\overline{\C}$ are (like $\mathcal{M}$) partially ordered rings with zero-divisors and we have
$\overline{\R}\subset\overline{\C}$. 
\end{remark}

	Due to the factorization we have: $v(a)=\infty$ \ifff $a=0$ in $\overline{\C}$. Consequently,
$\overline{\C}$ is a metric space (in contrast to
$\mathcal{M}$ is is a pseudometric space). The next two results follow immediately from the corresponding
properties of
$\mathcal{M}$ (Section~\ref{S: Functions of Moderate Growth}):

\begin{theorem}[Metric]\label{Metric} Colombeau's ring
$(\overline{\C}, d_v)$ is a ultrametric space, i.e. for every $a, b, c\in\overline{\C}$ we have:

	{\bf (i)}\quad $d_v(a, b)=0$ \ifff $a=b$.

	{\bf (ii)}\; $d_v(a, b) = d_v(b, a)$.

	{\bf (iii)} $d_v(a, b) \leq \max\{d_v(a, c), d_v(c, g)\}\leq d_v(a, c) + d_v(c, b)$. 
\end{theorem}

\begin{theorem}[Convergence in $\overline{\C}$]\label{T: Convergence in Cbar}  Let $(c_n)$ be
a sequence in
$\overline{\C}$ and
$c\in\overline{\C}$. Then:\newline
	\indent{\bf (i)} $\lim_{n\to\infty}c_n = 0$ in $(\overline{\C}, d_v)$ \ifff\,
$\lim_{n\to\infty}v(c_n) = \infty$ in
$\R\cup\{\infty\}$. {\bf (ii)} $\lim_{n\to\infty}c_n = c\not=0$ in $(\overline{\C}, d_v)$
implies that $v(c_n) = v(c)$ for all sufficiently large $n$.
\end{theorem}

		 The concept of a $v$-asymptotic series in $\mathcal{M}$ can be adapted to $\overline{\C}$. In what
follows $\lambda$ is the scale of $\overline{\C}$ (Definition~\ref{D: Colombeau's Ring}).

\begin{definition}[$v$-Asymptotic Series in $\overline{\mathbb{C}}$]

	{\bf (i)} A series
$\sum_{n=0}^\infty\frac{a_n}{\lambda^{r_n}}$ in $\overline{\C}$ is called {\bf $v$-asymptotic} if:

	\noindent (a) The set
of the coefficients $\{a_n\mid n\in\N\}$ is linearly $v$-independent over \C. 

	\noindent (b) The restriction of
the sequence $(r_n)$ on the set $\{n\mid a_n\not=0\}$ is a strictly increasing (finite or infinite)
sequence in \R.

	If $\mathbb{K}$ is a subset of $\mathcal{C}_v(\overline{\C})\cup\{0\}$, then
$\mathbb{K}(1/\lambda^\omega)$ denotes the set of all $v$-asymptotic series
$\sum_{n=0}^\infty\frac{a_n}{\lambda^{r_n}}$ in $\overline{\C}$ with coefficients $a_n$ in $\mathbb{K}$.

	{\bf (ii)} We say that $\sum_{n=0}^\infty\frac{a_n}{\lambda^{r_n}}$ is a {\bf $v$-asymptotic expansion}
of a generalized  number $a\in\overline{\C}$ (or, that $a$ is an {\bf asymptotic sum} of
$\sum_{n=0}^\infty\frac{a_n}{\lambda^{r_n}}$), in symbol
$a\leadsto\sum_{n=0}^\infty\frac{a_n}{\lambda^{r_n}}$, if\;
$\sum_{n=0}^\infty\frac{a_n}{\lambda^{r_n}}$ is a $v$-asymptotic series in $\overline{\C}$ and for every
$n\geq0,\; a_n\not=0$, we have $v\left(a-\sum_{k=0}^n\frac{a_k}{\lambda^{r_n}}\right)>r_n$.
\end{definition}
\begin{theorem}[Convergent Series in $\overline{\C}$]\label{C: Convergent Series in Cbar}
Let
$f\in\mathcal{M}$ and let
$\sum_{n=0}^\infty\,\frac{a_n}{x^{r_n}}$ be a $v$-asymptotic series in
$\overline{\mathbb{C}}$ such that $\lim_{n\to\infty}r_n=\infty$. Then the following are equivalent:
	
	{\bf (i)} $a\leadsto\sum_{n=0}^\infty\frac{a_n}{\lambda^{r_n}}$ in $\overline{\mathbb{C}}$.

	{\bf (ii)} $a=\sum_{n=0}^\infty\frac{a_n}{\lambda^{r_n}}$ in the sense that $a=\lim_{n\to\infty}\;
\sum_{k=0}^n\, \frac{a_k}{x^{r_k}}$ in $(\overline{\mathbb{C}}, d_v)$.
\end{theorem}
\Proof The result follows immediately from 
part~(i) of Theorem~\ref{T: Convergence in Cbar}. $\blacktriangle$
\begin{theorem}[$v$-Asymptotic Expansion in
$\overline{\C}$]\label{T: Existence and Uniqueness in Cbar} Every Colombeau's generalized number has a
$v$-asymptotic expansion in
$\overline{\mathbb{C}}$. More precisely, let
$C$ be a vector subspace of $\overline{\mathbb{C}}$ (over $\mathbb{C}$) which is 
maximal in $\mathcal{C}_v(\overline{\mathbb{C}})$. Then every
$a\in\overline{\mathbb{C}}$ has a unique
$v$-asymptotic expansion:
\begin{equation}\label{E: AsyNums}
a\leadsto\sum_{n=0}^\infty\frac{a_n}{\lambda^{r_n}},
\end{equation}
with coefficients $a_n$ in $C$.
\end{theorem}
\Proof We have $a=q(f)$ for some $f\in\mathcal{M}$. We observe that $\mathcal{C}=q^{-1}[C]$ is a
vector subspace of $\mathcal{M}$ which is maximal in $\mathcal{C}_v$ (Section~\ref{S: 
Pseudostandard Part}). It follows that
$f$ has a unique
$v$-asymptotic expansion $f(x)\leadsto\sum_{n=0}^\infty\frac{\varphi_n(x)}{x^{r_n}}$ in
$\mathcal{C}(1/x^\omega)$, by Theorem~\ref{T: Existence and Uniqueness} (Section~\ref{S: Existence and
Uniqueness Result}). Thus (\ref{E: AsyNums}) follows for $a_n=q(\varphi_n)$.
$\blacktriangle$

	If $\Omega$ is an open subset of $\R^n$, then we denote
\[
\mu_v(\Omega)=\{\omega+d\omega\mid \omega\in\Omega,\; d\omega\in\overline{\R}\,^n,\;
v(||d\omega||)>0 \}.
\] 
In what follows $\mathcal{G}(\Omega)$ denotes the algebra of Colombeau's generalized functions on
$\Omega$. We recall that $\mathcal{G}(\Omega)$ contains a canonical copy
of the space of Schwartz's distributions, in symbol,
$\mathcal{D}^\prime(\Omega)\subset\mathcal{G}(\Omega)$ (J.F. Colombeau~\cite{jfCol84}).
\begin{corollary} Let $f\in\mathcal{G}(\Omega)$ be a Colombeau's generalized function (in particular, a
Schwartz distribution) and let
$f:\mu_v(\Omega)\to\overline{\C}$ be its graph. Then there exist \newline $f_n:\mu_v(\Omega)\to
\overline{\C},\; f_n\not=0$, and $r_n:\mu_v(\Omega)\to\R$,  such that for every
$\xi\in\mu_v(\Omega)$: 

	(a) $\sum_{n=0}^\infty\frac{f_n(\xi)}{\lambda^{r_n(\xi)}}$ is a $v$-asymptotic series in
$\overline{\C}$ (which, among other things, implies that $v(f_n(\xi))=0$); 

	(b)
$f(\xi)\leadsto\sum_{n=0}^\infty\frac{f_n(\xi)}{\lambda^{r_n(\xi)}}$ in $\overline{\C}$.
\end{corollary}
\Proof The result follows immediately from Theorem~\ref{T: Existence and Uniqueness in Cbar} and the
identification between the Colombeau generalized functions with theirs graphs (M. Kunzinger and M.
Oberguggenberger~\cite{KunzOber99}).
$\blacktriangle$
%_____________________________________________________

}\end{document}